# Uniform Asymptotics for the Incomplete Gamma Functions Starting From Negative Values of the Parameters

N. M. Temme

CWI, P.O. Box 94079, 1090 GB Amsterdam, The Netherlands
*e-mail:* `nicot@cwi.nl`

**Abstract**

We consider the asymptotic behavior of the incomplete gamma functions $\gamma(-a, -z)$ and $\Gamma(-a, -z)$ as $a \to \infty$. Uniform expansions are needed to describe the transition area $z \sim a$, in which case error functions are used as main approximants. We use integral representations of the incomplete gamma functions and derive a uniform expansion by applying techniques used for the existing uniform expansions for $\gamma(a, z)$ and $\Gamma(a, z)$. The result is compared with Olver's uniform expansion for the generalized exponential integral. A numerical verification of the expansion is given in a final section.

## 1. Introduction

The incomplete gamma functions are defined by the integrals

$$\gamma(a,z) = \int_0^z t^{a-1} e^{-t}\, dt, \quad \Gamma(a,z) = \int_z^\infty t^{a-1} e^{-t}\, dt, \qquad (1.1)$$

where $a$ and $z$ are complex parameters and $t^a$ takes its principal value. For $\gamma(a,z)$ we need the condition $\Re a > 0$, for $\Gamma(a,z)$ we assume that $|\arg z| < \pi$. Analytic continuation can be based on these integrals, or on series representations of $\gamma(a,z)$. We have $\gamma(a,z) + \Gamma(a,z) = \Gamma(a)$.

Another important function is defined by

$$\gamma^*(a,z) = \frac{z^{-a}}{\Gamma(a)}\,\gamma(a,z) = \frac{1}{\Gamma(a)} \int_0^1 u^{a-1} e^{-zu}\, du. \qquad (1.2)$$

This function is a single-valued entire function of both $a$ and $z$, and is real for positive and negative values of $a$ and $z$. For $\Gamma(a,z)$ we have the additional integral representation

$$\Gamma(a,z) = \frac{e^{-z}}{\Gamma(1-a)} \int_0^\infty \frac{e^{-zt}\,t^{-a}}{t+1}\, dt, \quad \Re a < 1, \quad \Re z > 0, \qquad (1.3)$$

which can be verified by differentiating the right-hand side with respect to $z$.

For other information on the incomplete gamma functions we refer to Chapter IX of the BATEMAN MANUSCRIPT PROJECT (1953) or Chapter 11 of TEMME (1996).

The purpose of the paper is to derive new uniform asymptotic expansions for the functions $\gamma(-a,-z), \Gamma(-a,-z)$ and $\gamma^*(-a,-z)$. The uniform expansions for $\gamma(a,z)$ and $\Gamma(a,z)$ of our earlier papers, which will be summarized in the next section, are not valid for negative values of $a$.

Recently there is much interest in asymptotic properties of incomplete gamma functions. The function $\Gamma(a,z)$, or the related exponential integral

$$E_p(z) = z^{p-1} \int_z^\infty \frac{e^{-t}}{t^p} \, dt = z^{p-1} \Gamma(1-p, z), \tag{1.4}$$

plays an important role in Berry's smooth interpretation of the Stokes phenomena for certain integrals and special functions; see BERRY (1989). OLVER (1991A) investigated $E_p(z)$ in particular at the Stokes lines $\arg z = \pm \pi$ and used the results in OLVER (1991B); see also OLVER (1994). We summarize Olver's result in the next section. In DUNSTER (1994A) a new expansion for the function $E_p(z)$ is given for complex values of $z$ and large positive values of $p$ and in DUNSTER (1994B) error bounds are given, in particular for the Stokes' smoothing approximations. In PARIS (1994) the uniform asymptotic expansions of the next section are used for complex values of the parameters in a new algorithm for computing the Riemann zeta function on the critical line.

The computational problem of the incomplete gamma functions for complex values of $a$ and $z$ is not well solved in the software literature. In particular when the complex parameters have large negative real parts existing computer programs may give really false answers. We expect that the new expansions for $\gamma(-a,-z), \Gamma(-a,-z), \gamma^*(-a,-z)$ and Olver's expansion for $E_p(-z)$ will be of value to solve this problem.

## 2. Uniform expansions for incomplete gamma functions.

We summarize the known uniform expansions for the incomplete gamma functions; see TEMME (1975), (1979) and more recently OLVER (1991A), (1994).

*2.1. Uniform expansions for P and Q.*

Let $\eta$ be the real number defined by

$$\tfrac{1}{2}\eta^2 = \lambda - 1 - \ln \lambda, \quad \lambda > 0, \quad \text{sign}(\eta) = \text{sign}(\lambda - 1). \tag{2.1}$$

Extend the relation between $\lambda$ and $\eta$ to complex values by analytic continuation. We have

$$\eta \sim (\lambda - 1)\left[1 - \tfrac{1}{3}(\lambda - 1) + \ldots\right], \quad \lambda \to 1.$$

Then we write, with $\lambda = z/a$,

$$\begin{aligned} Q(a,z) &= \frac{\Gamma(a,z)}{\Gamma(a)} = \tfrac{1}{2}\,\text{erfc}(\eta\sqrt{a/2}\,) + R_a(\eta), \\ P(a,z) &= \frac{\gamma(a,z)}{\Gamma(a)} = \tfrac{1}{2}\,\text{erfc}(-\eta\sqrt{a/2}\,) - R_a(\eta). \end{aligned} \tag{2.2}$$

The error functions are defined by

$$\operatorname{erf} z = \frac{2}{\sqrt{\pi}} \int_0^z e^{-t^2} \, dt, \quad \operatorname{erfc} z = 1 - \operatorname{erf} z = \frac{2}{\sqrt{\pi}} \int_z^\infty e^{-t^2} \, dt \qquad (2.3)$$

The error functions are the dominant terms in (2.2) as $a$ tends to infinity, and they describe the transition at $a = z$.

We have

$$R_a(\eta) = \frac{e^{-\frac{1}{2}a\eta^2}}{\sqrt{2\pi a}} S_a(\eta), \quad S_a(\eta) \sim \sum_{n=0}^\infty \frac{C_n(\eta)}{a^n}, \qquad (2.4)$$

as $a \to \infty$, uniformly with respect to $|\arg a| \leq \pi - \delta_1$, and $\lambda, |\arg \lambda| \leq 2\pi - \delta_2$, where $\delta_1, \delta_2$ are arbitrarily small positive constants. More information on the coefficients $C_n(\eta)$ is given in the next subsection.

In this paper we derive similar new expansions for $\gamma(a, z)$ and $\Gamma(a, z)$, where we concentrate on negative values of $a$ and $z$. Because simple reflection formulas are missing that relate $\gamma(-a, z)$ with $\gamma(a, z)$, etc., we cannot simply transform the above expansions to this case.

*2.2. Further details on the previous results.*

Because the new expansions (see §3) are closely related with the one in (2.4), we give some details on computing the coefficients and on the mapping $\lambda \to \eta(\lambda)$ defined in (2.1).

The expansion in (2.4) for $S_a(\eta)$ can be obtained by differentiating one of the equations in (2.2) with respect to $\eta$, which gives

$$\frac{d}{d\eta} S_a(\eta) - a\eta S_a(\eta) = a \left[ 1 - f(\eta)/\Gamma^*(a) \right],$$

where

$$f(\eta) = \frac{1}{\lambda} \frac{d\lambda}{d\eta} = \frac{\eta}{\lambda - 1}, \quad \Gamma^*(a) = \sqrt{a/(2\pi)} \, e^a \, a^{-a} \, \Gamma(a). \qquad (2.5)$$

Substituting the asymptotic expansion for $S_a(\eta)$ one finds for the coefficients the relations

$$C_0(\eta) = \frac{1}{\lambda - 1} - \frac{1}{\eta}, \quad \eta C_n(\eta) = \frac{d}{d\eta} C_{n-1}(\eta) + \gamma_n f(\eta), \quad n \geq 1, \qquad (2.6)$$

where $\gamma_n$ are the coefficients in the reciprocal gamma function expansion

$$\frac{1}{\Gamma^*(a)} \sim \sum_{n=0}^\infty \frac{\gamma_n}{a^n} \quad a \to \infty, \quad |\arg a| < \pi, \qquad (2.7)$$

of which the first few are given by

$$\gamma_0 = 1, \quad \gamma_1 = -\frac{1}{12}, \quad \gamma_2 = \frac{1}{288}, \quad \gamma_3 = \frac{139}{51840}, \quad \gamma_4 = -\frac{571}{2488320}.$$

Next we summarize the properties of the mapping $\lambda \to \eta(\lambda)$ defined in (2.1). More details are given in TEMME (1979). Write $\eta = \alpha + i\beta$. The mapping is one-to-one for $-2\pi < \arg \lambda < 2\pi$, $\lambda \neq 0$. The corresponding $\eta$–domain is given in Figure 1. The mapping is singular at $\lambda = \exp(\pm 2\pi i)$ with corresponding points $B^{\pm} = \sqrt{2\pi}\,(-1 \pm i)$. The positive real $\lambda$-axis (with $\arg \lambda = \pm 2\pi$) is mapped to branches of the hyperbolas $\alpha\beta = \pm 2\pi$, $\alpha < -\sqrt{2\pi}$. When we cut the $\eta$-plane along the two thick parts of the hyperbolas in Figure 1 we have a one-to-one interpretation of the mapping in (2.1), and the function $\eta(\lambda)$ is singular at $\eta = \sqrt{2\pi}\,(-1 \pm i)$. The parabola shaped curve in the $\lambda$-plane is the set of points $\Im(\lambda - \ln \lambda - 1) = 0$, and given by the equation $\rho(\theta) = \theta/\sin\theta$, $-\pi < \theta < \pi$, $\rho(0) = 1$, where $\rho, \theta$ are the polar coordinates of $\lambda$: $\lambda = \rho \exp(i\theta)$. The corresponding set in the $\eta$-plane is the imaginary axis.

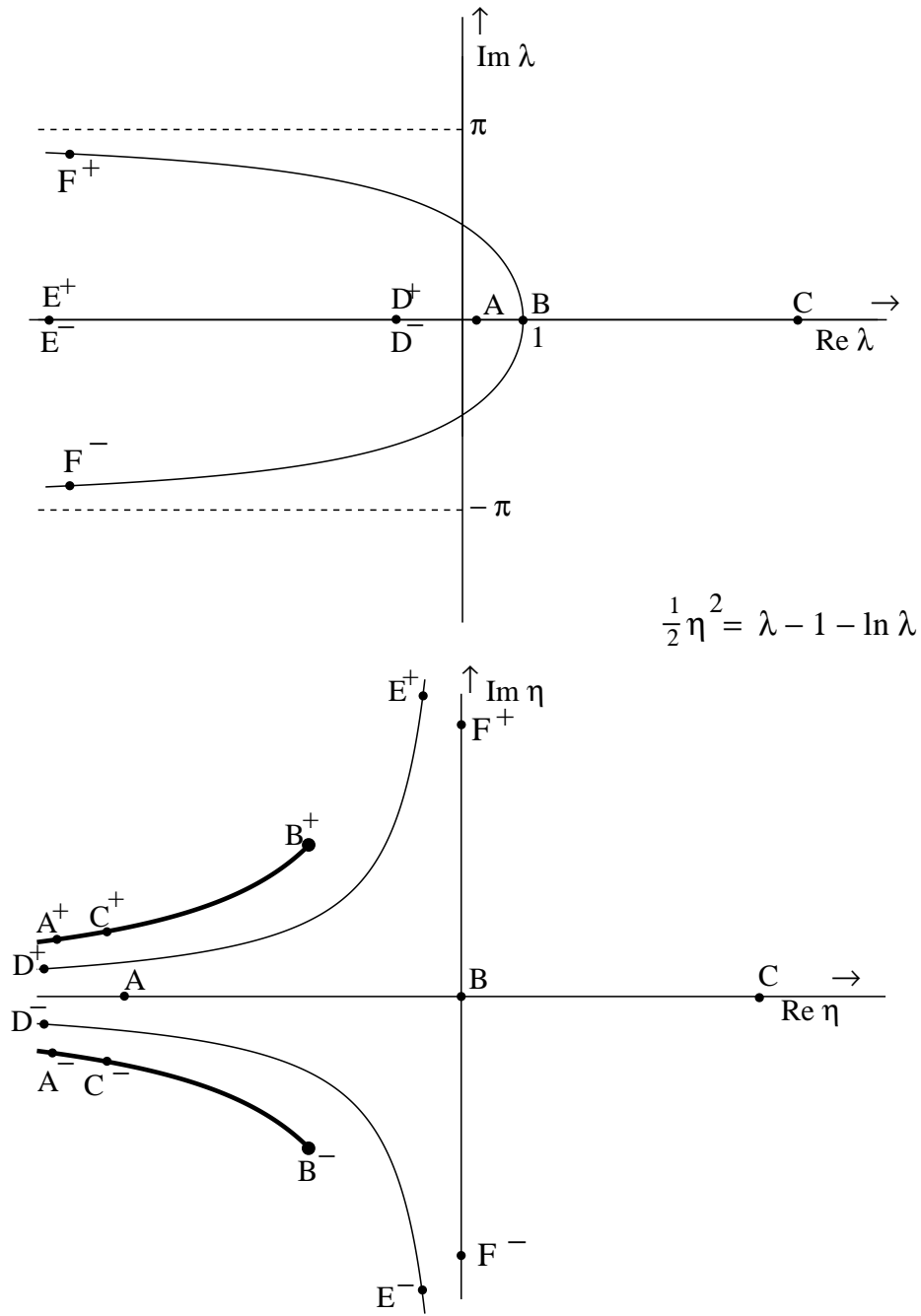

**Figure 1.** Corresponding points in the $\lambda-$ and $\eta-$planes for the mapping $\lambda \to \eta(\lambda)$ defined in (2.1). The points $B^{\pm}$ are singular points in the $\eta-$plane. In the $\lambda-$plane they correspond to $\exp(\pm 2\pi i)$. The $\eta-$points $D^{\pm}, E^{\pm}$ correspond to $\lambda-$points $D^{\pm}, E^{\pm}$ having phases $\pm \pi i$. The thick branches of the hyperbolas are branch cuts. The parabola shaped curve in the $\lambda-$plane corresponds to the imaginary axis in the $\eta-$plane.

*2.3. Olver's expansion.*

OLVER (1991A) investigated the generalized exponential integral $E_p(z)$ defined in (1.4). Let

$$F_p(z) = \frac{\Gamma(p)}{2\pi} \frac{E_p(z)}{z^{p-1}}$$

and $z = \rho e^{i\theta}$, $\alpha = n - \rho + p$ with $\rho$ a large parameter, $p$ fixed. Then

$$F_{n+p}(z) \sim (-1)^n i e^{-\rho \pi i} \left[ \tfrac{1}{2}\mathrm{erfc}\left\{ c(\theta)\sqrt{\tfrac{1}{2}\rho} \right\} - i \frac{e^{-\frac{1}{2}\rho\{c(\theta)\}^2}}{\sqrt{2\pi\rho}} \sum_{s=0}^{\infty} \left(\tfrac{1}{2}\right)_s g_{2s}(\theta, \alpha) \left(\tfrac{2}{\rho}\right)^s \right], \quad (2.8)$$

uniformly with respect to $\theta \in [-\pi + \delta, 3\pi - \delta]$ and bounded values of $|\alpha|$; $\delta$ denotes an arbitrarily small positive constant. Furthermore,

$$c(\theta) = \sqrt{2\left\{ e^{i\theta} + i(\theta - \pi) + 1 \right\}},$$

with the choice of branch of the square root that implies $c(\theta) \sim (\pi - \theta)$ as $\theta \to \pi$; the coefficients $g_{2s}(\theta, \alpha)$ are continuous functions of $\theta$ and $\alpha$. A similar expansion for $F_{n+p}(z)$ is given when $\theta \in [-3\pi + \delta, \pi - \delta]$. As Olver remarks, this expansion quantifies the Stokes phenomenon, that is, the rapid but smooth change in form of other expansions as $\theta$ passes through the common interval of validity of the other expansions.

In the present paper we give an expansion that is also valid in a domain around the negative $z$-axis and that also contains the error function. Our expansion is quite different, however, and is strongly related to the uniform expansion given in (2.1)–(2.4).

## 3. Starting from negative values.

We start with (1.3) and replace $a$ with $-a$:

$$\Gamma(-a, z) = \frac{e^{-z}}{\Gamma(1+a)} \int_0^{\infty} \frac{e^{-zt}\, t^a}{t+1}\, dt, \quad \Re a > -1, \quad \Re z > 0. \tag{3.1}$$

By turning the path of integration and invoking the principle of analytic continuation we can enlarge the domain of $z$. For example, when we consider the path in (3.1) along the positive imaginary axis, the integral is defined for $-\pi < \arg z < 0$, and in the overlapping domain $-\pi/2 < \arg z < 0$ its value is the same as in (3.1). We turn the path from the positive imaginary axis $\arg t = \pi/2$ to the negative axis $\arg t = \pi$, avoiding the pole at $t = -1$ by using a small semi-circle. The integral is then defined for $-\frac{3}{2}\pi < \arg z < -\frac{1}{2}\pi$. We change the variable of integration and the result is:

$$\Gamma(-a, ze^{-i\pi}) = \frac{e^z e^{i\pi a}}{\Gamma(1+a)} \int_0^\infty \frac{e^{-zt} t^a}{t-1} dt, \quad \Re a > -1, \quad -\frac{1}{2}\pi < \arg z < \frac{1}{2}\pi \qquad (3.2)$$

and we avoid the pole by integrating under the pole at $t = 1$. By turning the path in (3.1) clockwise we obtain:

$$\Gamma(-a, ze^{+i\pi}) = \frac{e^z e^{-i\pi a}}{\Gamma(1+a)} \int_0^\infty \frac{e^{-zt} t^a}{t-1} dt, \quad \Re a > -1, \quad -\frac{1}{2}\pi < \arg z < \frac{1}{2}\pi \qquad (3.3)$$

and we avoid the pole by integrating above the pole at $t = 1$. An easy consequence is

$$e^{i\pi a}\Gamma(-a, ze^{+i\pi}) - e^{-i\pi a}\Gamma(-a, ze^{-i\pi}) = \frac{-2\pi i}{\Gamma(a+1)}, \qquad (3.4)$$

which follows from computing the residue of the integral over the full circle around the pole at $t = 1$. It also follows from the fact that $\Gamma(a, z) = \Gamma(a)[1 - z^a \gamma^*(a, z)]$ and $\gamma^*(a, z)$ is entire.

We proceed with (3.3) and assume, for the time being, that $a$ and $z$ are positive. We change the variable of integration $t \to at/z$, and obtain

$$\Gamma\left(-a, ze^{\pi i}\right) = \frac{e^z \lambda^{-a} e^{-ia\pi}}{\Gamma(a+1)} \int_0^\infty e^{-at} t^a \frac{dt}{t-\lambda}$$

$$= \frac{e^{\frac{1}{2}a\eta^2 - ia\pi}}{\Gamma(a+1)} \int_{-\infty}^\infty e^{-\frac{1}{2}a\zeta^2} g(\zeta) \frac{d\zeta}{\zeta - \eta},$$

where $\eta = \eta(\lambda)$ is defined in (2.1), $\lambda = z/a$, $\zeta = \eta(t)$, that is,

$$\frac{1}{2}\zeta^2 = t - \ln t - 1, \quad \mathrm{sign}\zeta = \mathrm{sign}(t-1),$$

and

$$g(\zeta) = \frac{dt}{d\zeta} \frac{\zeta - \eta}{t - \lambda} = \frac{\zeta t}{t-1} \frac{\zeta - \eta}{t - \lambda}. \qquad (3.5)$$

In the $\zeta$-integral the path passes above the pole at $\zeta = \eta$.

When $z \sim a$, that is, $\lambda \sim 1$, the pole is near the saddle point $\zeta = 0$, and for large values of $a$ we need an error function to describe the asymptotic behavior. We can split off the pole by writing $g(\zeta) = [g(\zeta) - g(\eta)] + g(\eta)$, where, as is easily verified, $g(\eta) = 1$. Using the representation

$$\tfrac{1}{2}\mathrm{erfc}\, iz = -\frac{e^{z^2}}{2\pi i} \int_{-\infty}^\infty \frac{e^{-t^2}}{t-z} dt \quad z \in \mathbb{C},$$

where the path passes above the pole at $t = z$, we obtain (introducing a suitable normalization)

$$-\Gamma(a+1)\frac{e^{\pi i a}}{2\pi i} \Gamma\left(-a, ze^{\pi i}\right) = \tfrac{1}{2}\mathrm{erfc}\left(i\eta\sqrt{\tfrac{a}{2}}\right) - i\frac{e^{\frac{1}{2}a\eta^2}}{\sqrt{2\pi a}} T_a(\eta), \qquad (3.6)$$

and, for instance by using the relation in (3.4),

$$\Gamma(a+1)\frac{e^{-\pi i a}}{2\pi i}\Gamma\left(-a, ze^{-\pi i}\right) = \tfrac{1}{2}\mathrm{erfc}\left(-i\eta\sqrt{\tfrac{a}{2}}\right) + i\frac{e^{\frac{1}{2}a\eta^2}}{\sqrt{2\pi a}}T_a(\eta), \qquad (3.7)$$

where

$$T_a(\eta) = -\sqrt{\frac{a}{2\pi}}\int_{-\infty}^{\infty} e^{-\frac{1}{2}a\zeta^2} h(\zeta)\,d\zeta, \quad h(\zeta) = \frac{g(\zeta) - g(\eta)}{\zeta - \eta}. \qquad (3.8)$$

Standard methods for integrals can be used now (see WONG (1989)) to obtain an asymptotic expansion of $T_a(\eta)$ in negative powers of $a$. It is easier, however, to use a differential equation satisfied by $T_a(\eta)$. In this way we can identify the coefficients of the expansion with those in (2.4).

Differentiating (3.6) with respect to $\eta$, and using (1.1), (2.1) and (2.3) we obtain the differential equation

$$\frac{d}{d\eta}T_a(\eta) + a\eta T_a(\eta) = a\left[f(\eta)\Gamma^*(a) - 1\right], \qquad (3.9)$$

where $f(\eta)$ and $\Gamma^*(a)$ are given in (2.5). As in (2.7) we have (with the same coefficients $\gamma_n$ indeed)

$$\Gamma^*(a) \sim \sum_{n=0}^{\infty}(-1)^n\frac{\gamma_n}{a^n}, \quad a\to\infty, \quad |\arg a| < \pi.$$

Substituting this and the expansion

$$T_a(\eta) \sim \sum_{n=0}^{\infty}(-1)^n\frac{C_n(\eta)}{a^n}, \qquad (3.10)$$

in (3.9) we find for $C_n(\eta)$ the same recursion as in (2.6), and we conclude that in (2.4) and (3.10) the same coefficients occur.

*3.1. The domain of validity of expansion (3.10).*

We return to the case that $a$, $z$ and $\lambda = z/a$ are complex parameters and claim that expansion (3.10) holds as $a \to \infty$, uniformly with respect to $|\arg a| \leq \pi - \delta_1$ and $|\arg \lambda| \leq 2\pi - \delta_2$, with $\delta_1, \delta_2$ arbitrarily small positive constants.

To verify this we observe that $\eta(\lambda)$ is analytic and univalent in the sector $|\arg \lambda| < 2\pi$. This follows from §2.2. The function $\zeta(t)$ is in fact $\eta(t)$ and the quantity $t = t(\zeta)$ occurring in $g$ in (3.5) is singular at the points $B^{\pm}$ shown in the lower part of Figure 1. The functions $g$ and $h$ (given in (3.8)) are singular in the same points. It follows that $h(\zeta)$ is analytic in the disk $|\zeta| < 2\sqrt{\pi}$ and in the sectors $|\Im\zeta| < |\Re\zeta|$. When $a$ is complex we can turn the path of integration in (3.8) inside the sectors $|\Im\zeta| < |\Re\zeta|$ and $\Re a\zeta^2$ can be kept positive on the path of integration as long as $|\arg a| < \pi$. This gives the domains for $a$ and $\lambda$.

We remark that the function $h(\zeta)$ is bounded at infinity inside the sectors $|\Im\zeta| < |\Re\zeta|$. This property is not needed to verify the above result, however.

*3.2. A real expansion for $\gamma^*(-a,-z)$.*

It is of interest to have a result for the function $\gamma^*(a,z)$ (see (1.2)) for negative values of the parameters. In that case $\gamma^*(a,z)$ is real and we have a real asymptotic representation:

$$\gamma^*(-a,-z) = z^a \left\{ \cos \pi a - \sqrt{\tfrac{2a}{\pi}}\, e^{\frac{1}{2}a\eta^2}\, \sin \pi a \, \left[ \sqrt{\tfrac{2}{a}}\, F\left(\eta\sqrt{\tfrac{a}{2}}\right) + \tfrac{1}{a}T_a(\eta) \right] \right\}, \qquad (3.11)$$

where

$$F(z) = e^{-z^2} \int_0^z e^{t^2}\, dt$$

which is Dawson's integral. The relation with the error function is:

$$F(z) = -\tfrac{1}{2}i\sqrt{\pi}\, e^{-z^2}\, \mathrm{erf}\, iz.$$

Representation (3.11) can be verified by using (3.3) and the relations

$$\Gamma(-a, ze^{+i\pi})\Gamma(-a)\left[1 - z^{-a}e^{-i\pi a}\gamma^*(-a,-z)\right], \quad \Gamma(a+1)\Gamma(-a) = -\frac{\pi}{\sin \pi a}.$$

In (3.11) we see that the oscillatory behavior of $\gamma^*(-a,-z)$ is described by two terms, one with $\cos \pi a$ and one with $\exp(\tfrac{1}{2}a\eta^2)\sin \pi a$ (the factor containing Dawson's integral is slowly varying when the parameters $a$ and $z$ are positive). This complicated oscillatory behavior is one of the problems in writing reliable software for the functions $\gamma(a,z)$ and $\Gamma(a,z)$ when the parameters have large negative real parts. Dawson's integral becomes dominant when we consider complex values of the parameters. The function $F(z)$ is an entire odd function and has for $\Re z \geq 0$ the asymptotic behavior

$$F(z) \sim \begin{cases} \frac{1}{2z} & \text{if } \arg z < \tfrac{1}{4}\pi, \\ \mathrm{sign}(\Im z)\tfrac{1}{2}i\sqrt{\pi}\, e^{-z^2} & \text{elsewhere.} \end{cases} \quad |z| \to \infty,$$

*3.3. Comparing Olver's expansion with our expansion.*

In Olver's expansion (2.9) the modulus of $z$ is the large parameter, whereas in our expansion (3.10) $a$ is the large parameter. In addition, in Olver's expansion the transition occurs at $\theta = \pi$, whereas in our expansion the transition takes place at $\eta = 0$, that is, at $\lambda = 1$, with $\lambda = z/a$. We have $\eta \sim \lambda - 1$, as $\lambda \to 1$. Hence, when $z = a + i\varepsilon$, then $\eta \sim i\varepsilon/a$. When $a$ is positive and $\varepsilon$ increases from negative to positive values, $ze^{\pi i}$ crosses the negative axis with increasing argument, and $\tfrac{1}{2}\mathrm{erfc}(i\eta\sqrt{a/2})$ of (3.6) changes rapidly from 0 to 1.

Olver's expansion has more parameters than ours because his analysis started by considering $F_{n+p}(z)$ as a remainder in the expansion for $F_p(z)$. However, one can take $n = 0, p = a+1$. Then Olver's notation agrees better with our notation. In his result the parameter $\alpha = \rho - p$ is assumed to be bounded, and a restriction like this is not present in our method.

**Remark.** We verify what happens when the parameters $a, z$ in (2.2) are taken with negative signs. First let $\Im a, \Im z$ be positive, and replace in $Q(a, z)$ the parameters $a, z$ with $ae^{-\pi i}, ze^{-\pi i}$, respectively. The quantity $\eta$ defined in (2.1) does not change by this operation, whereas the expansion for $S_a(\eta)$ becomes the expansion for $T_a(\eta)$ given in (3.10). When we formally write $S_{-a}(\eta) = T_a(\eta)$, which certainly is not true, it follows that the right-hand side of the first line in (2.2) formally becomes

$$\tfrac{1}{2}\mathrm{erfc}\left(-i\eta\sqrt{a/2}\right) + i\frac{e^{\frac{1}{2}a\eta^2}}{\sqrt{2\pi a}}T_a(\eta),$$

which is the right-hand side of (3.7). But

$$Q(-a, ze^{-\pi i}) = \frac{\Gamma\left(-a, ze^{-\pi i}\right)}{\Gamma(-a)} = -\frac{1}{\pi}\Gamma(a+1)\sin\pi a\,\Gamma(-a, ze^{-\pi i}),$$

which means that the left-hand side of (3.7) equals

$$\Gamma(a+1)\frac{e^{-\pi i a}}{2\pi i}\Gamma\left(-a, ze^{-\pi i}\right) = \frac{1}{1 - e^{2\pi i a}}\,Q\left(-a, ze^{-\pi i}\right).$$

Hence, by proceeding formally from the relation for $Q$ in formula (2.2), we miss the factor $1/[1 - \exp(2\pi i a)]$, this factor being negligible when $\Im a$ is positive and large, as we assumed here. A similar conclusion holds when $\Im a, \Im z$ be negative and we replace $a, z$ with $ae^{+\pi i}, ze^{+\pi i}$.

## 4. Numerical verification of the results.

To verify numerically the uniform expansion we have used the real representation (3.11) for $\gamma^*(-a, -z)$ with large positive values of $a$ and $z$. We have used the recursion formula

$$\gamma^*(-a-1, -z) + z\gamma^*(-a, -z) = -\frac{1}{\pi}\sin\pi a\,e^z\Gamma(a+1) \tag{4.1}$$

for testing the results. This relation easily follows from the well-known recursion $a\gamma(a, z) - \gamma(a+1, z) = z^a e^{-z}$ and the relation between $\gamma(a, c)$ and $\gamma^*(a, z)$ given in (1.2). To avoid overflow and strong oscillations we have normalized the function $\gamma^*(a, z)$ by using a quantity $\widetilde{\gamma}_a(z)$ defined by

$$\gamma^*(-a, -z) = z^a\,\cos\pi a + \sin\pi a\,\Gamma(a)\,e^z\,\widetilde{\gamma}_a(z). \tag{4.2}$$

The recursion reads for this new function reads

$$-\widetilde{\gamma}_{a+1}(z) + \lambda\,\widetilde{\gamma}_a(z) + \frac{1}{\pi} = 0, \quad \lambda = \frac{z}{a}. \tag{4.3}$$

Large quantities in the right-hand side of (3.11) are $z^a$ and $z^a\exp(\tfrac{1}{2}a\eta^2)$, and when dividing (3.11) by $e^z\Gamma(a)$ we can use the relation

$$\frac{z^a\,e^{\frac{1}{2}a\eta^2}}{e^z\Gamma(a)} = \sqrt{\frac{a}{2\pi}}\,\frac{1}{\Gamma^*(a)},$$

where $\Gamma^*(a)$ is defined in (2.5). This gives

$$\widetilde{\gamma}_a(z) = -\frac{a}{\pi\,\Gamma^*(a)}\left[\sqrt{\tfrac{2}{a}}\,F\left(\eta\sqrt{\tfrac{a}{2}}\right) + \tfrac{1}{a}T_a(\eta)\right], \qquad (4.4)$$

Dawson's function $F$ and $\Gamma^*(a)$ are computed with extended precision (about 19 relevant digits). We have used expansion (3.10) in (4.4) and truncated the series after the term $C_6(\eta)/a^6$, with the expectation that the order of magnitude of the remainder is about $10^{-14}$ if $a \geq 100$.

We have used the representations of the coefficients $C_0, \ldots C_6$ as given in Section 3 of TEMME (1979). For $|\eta| \leq 1$ we have used Maclaurin expansions of the coefficients. An extensive set of Maclaurin coefficients is given in DIDONATO & MORRIS (1986). We have used enough Maclaurin coefficients in order to obtain 14 digits accuracy for $C_0$ on the interval $|\eta| \leq 1$, 12 digits for $C_1$, etc., with the intention to use the asymptotic expansion if $a \geq 100$.

In the numerical verification we fixed $a = 100$ and took $z$ in intervals around the point $a$, the transition point. In fact we used random numbers for $z$ in intervals $[(1-2^{-k})a, (1+2^{-k})a]$, for $k = 0, 1, \ldots 6$. All computed left-hand sides of the recursion relation (4.3) were (in modulus) less than $5.5 \times 10^{-18}$. Values of $\widetilde{\gamma}_a(z)$ and $\gamma^*(-a, -z)$ around the turning point $a = z = 100$ are given in Table 4.1, together with the error in the computation of $\widetilde{\gamma}_a(z)$.

| $z$ | $\widetilde{\gamma}_a(z)$ | $\gamma^*(-a,-z)$ | error |
|---|---|---|---|
| 90.0 | 2.464370784806670 | 1.20552423411674426e+196 | 5.3e–18 |
| 91.0 | 2.372768700058633 | 3.39249458130961502e+196 | 2.7e–19 |
| 92.0 | 2.245263163760808 | 9.45984032045937199e+196 | 5.4e–20 |
| 93.0 | 2.081986332010084 | 2.61247251775844538e+197 | 7.0e–19 |
| 94.0 | 1.884270784695821 | 7.14192476088520721e+197 | 5.4e–20 |
| 95.0 | 1.654640974379185 | 1.93183088802899244e+198 | 3.8e–19 |
| 96.0 | 1.396740233402602 | 5.16774449595400219e+198 | 1.4e–18 |
| 97.0 | 1.115197442856446 | 1.36641025407873601e+199 | 4.9e–19 |
| 98.0 | 0.815441855478639 | 3.56901459711775166e+199 | 1.1e–19 |
| 99.0 | 0.503478153282224 | 9.20223634675090428e+199 | 1.2e–18 |
| 100.0 | 0.185636311520584 | 2.34010604791689845e+200 | 7.0e–19 |
| 101.0 | −0.131687939165396 | 5.86259361067072689e+200 | 8.1e–19 |
| 102.0 | −0.442286347630879 | 1.44483553329971751e+201 | 9.5e–19 |
| 103.0 | −0.740370670716527 | 3.49587466488947244e+201 | 4.6e–19 |
| 104.0 | −1.020772405557316 | 8.28106485635440568e+201 | 8.7e–19 |
| 105.0 | −1.279100681803916 | 1.91262956371154387e+202 | 5.4e–20 |
| 106.0 | −1.511852726452274 | 4.28009100876763689e+202 | 1.2e–18 |
| 107.0 | −1.716474855286176 | 9.18429865111133800e+202 | 2.0e–18 |
| 108.0 | −1.891375273470379 | 1.85454711082400020e+203 | 4.9e–19 |
| 109.0 | −2.035892840388064 | 3.38697191526660474e+203 | 2.7e–19 |
| 110.0 | −2.150228198004453 | 5.01390135464872193e+203 | 1.0e–18 |

**Table 4.1.** Values of $\widetilde{\gamma}_a(z)$ of (4.2) and $\gamma^*(-a,-z)$ for $a = 100.25$ with $z$−values around the transition point $a = z$; the "error" is the computed left-hand side of (4.3). Notation: 1.20552423411674426e+196 means $1.20552423411674426 \times 10^{196}$, 5.3e–18 means $5.3 \times 10^{-18}$.